\documentclass[twoside,10.5pt]{article}
\usepackage{mathrsfs}
\usepackage{pifont}
\usepackage{amsmath}
\usepackage{amsthm}
\usepackage{txfonts}
 \usepackage{geometry}
\usepackage{latexsym}
\usepackage{amssymb}
\usepackage{graphicx}
\usepackage{geometry}
\usepackage{enumitem}%
\usepackage{ifthen}
\usepackage{fancyhdr}
\usepackage{afterpage}
\usepackage{xcolor}
\usepackage{ifthen}
\usepackage{titlesec}
\titleformat{\section}{\bfseries}{\thesection.}{0.3em}{}
\titlespacing{\section}{0pt}{0pt}{0pt}
\titleformat{\subsection}{\itshape}{\thesubsection}{0.3em}{}
\titlespacing{\subsection}{0pt}{0pt}{0pt}
\titleformat{\subsubsection}{\rmfamily}
{\thesubsubsection}{0.3em}{}
\titlespacing{\subsubsection}{0pt}{0pt}{0pt}
\newcommand{\eqdef}{:=} %
\newcommand{\sq}{\qquad $\square$}
\begin{document}

\begin{center}
{\LARGE{Monotonically Decreasing Sequence of Divergences}}\\[15pt]
Tomohiro Nishiyama\\
E-mail: htam0ybboh@gmail.com
\end{center}
\textbf{Abstract}

Divergences are quantities that measure discrepancy between two probability distributions and play an important role in various fields such as statistics and machine learning.
Divergences are non-negative and are equal to zero if and only if two distributions are the same. In addition, some important divergences such as the f-divergence have convexity, which we call ``convex divergence''. In this paper, we show new properties of the convex divergences by using integral and differential operators that we introduce. For the convex divergence, the result applied the integral or differential operator is also a divergence. In particular, the integral operator preserves convexity.
Furthermore, the results applied the integral operator multiple times constitute a monotonically decreasing sequence of the convex divergences.
We derive new sequences of the convex divergences that include the Kullback-Leibler divergence or the reverse Kullback-Leibler divergence from these properties. 


\textbf{1. Introduction}
Let $(\Sigma, \mathcal{M})$ be a measurable space (Billingsley, 2008) where $\Sigma$ denotes the sample space and $\mathcal{M}$ denotes the $\sigma$-algebra on $\Sigma$.
Let $\mathcal{P}$ denotes the set of probability distributions (measures) with common support on $(\Sigma, \mathcal{M})$.

Divergences are quantities that measure discrepancy between probability distributions $P,Q\in \mathcal{P}$, and they are defined as functions that satisfy the following properties (Amari {\&}  Cichocki, 2010).

Let $D: \mathcal{P} \times \mathcal{P}\rightarrow[0,\infty] $. 
For any $P,Q\in \mathcal{P}$,
\begin{align}
D(P\|Q)\geq 0 \nonumber \\ 
D(P\|Q)=0 \iff P=Q. \nonumber
\end{align}

The $f$-divergence (Ajjanagadde, Makur, Klusowski {\&} Xu, 2017; Csisz{\'a}r {\&} Shields, 2004; Sason {\&} Verdu, 2016) and the Bregman divergence (Bregman, 1967; Nishiyama, 2018) are well-known classes of divergences. 
For a strictly convex function $f: (0,\infty) \rightarrow \mathbb{R}$ and $f(1)=0$, the $f$-divergence is defined as 
\begin{align}
D_f(P\|Q)\eqdef\int qf\biggl(\frac{p}{q} \biggr)d\mu, \nonumber
\end{align}
where $p$ and $q$ are the Radon-Nikodym derivatives of $P$ and $Q$, respectively and $\mu$ denotes a positive measure (e.g. the Lebesgue measure or the counting measure).
The integral is defined on the support of $P$.
For a differentiable strictly convex function $F: (0,\infty) \rightarrow \mathbb{R}$, the Bregman divergence is defined as
\begin{align}
B_F(P\|Q)\eqdef \int \biggl(F(p)-F(q) - F'(q)(p-q) \biggr) d\mu, \nonumber
\end{align}
where $F'(x)$ denotes the derivative with respect to $x$.
For example, the Kullback-Leibler divergence (KL-divergence) $\mathrm{KL}(P\|Q)\eqdef \int p\log \frac{p}{q} d\mu$ (Kullback, 1997; Cover {\&} Thomas, 2012) is the most fundamental divergence and it belongs to both the $f$-divergence and the Bregman divergence, where $\log$ denotes the natural logarithm.
The $f$-divergence is convex in the first and the second argument and the Bregman divergence is convex in the first argument.
The convexity in the second argument is defined as follows.
For probability distributions $P, Q_1, Q_2\in\mathcal{P}$ and a parameter $\lambda\in[0,1]$, 
\begin{align}
\lambda D(P\|Q_1) + (1-\lambda)  D(P\|Q_2) \geq D(P\|\lambda Q_1 + (1-\lambda)Q_2). \nonumber
\end{align}
The convexity in the first argument is defined in the same way.
When divergences are convex in the first or the second argument, we call these divergence ``convex divergence''.
 
In this paper, we introduce integral and differential operators and show the following properties.

(a) The result applied the integral operator for a convex divergence is a convex divergence of which value is less than or equal to the original divergence.

(b) The convex divergences applied the integral operator multiple times constitute a monotonically decreasing sequence.

(c) The result applied the differential operator for a convex divergence is a divergence of which value is larger than or equal to the original divergence.

Furthermore, we show a specified examples applied these results to the KL-divergence and the reverse KL-divergence. 
Although the R{\'e}nyi divergences $D_\alpha(P\|Q)\eqdef\frac{1}{\alpha-1} \log \int p^\alpha q^{1-\alpha} d\mu$ constitute a divergence sequence that includes the KL-divergence, this sequence is a monotonically increasing (Van Erven {\&} Harremos, 2014).
In contrast, we derive new monotonically decreasing sequences of divergences including the KL-divergence or the reverse KL-divergence by using the polylogarighm (Wood, 1992).

\textbf{2. Definitions}

We show some definitions and notations in this paper.
Consider probability distributions $P,Q \in\mathcal{P}$ with common support.
Let $R(t)\in\mathcal{P}$ be a mixture distribution defined as $R(t)\eqdef (Q-P)t+P$ for $t\in[0,1]$.

If a divergence $D(P\|Q)$ is convex in the first argument, we call it ``left-convex divergence'' and if a divergence $D(P\|Q)$ is convex in the second argument, we call it  ``right-convex divergence''.

If a divergence $D(P\|Q)$ is left-convex, by putting $\hat{D}(P\|Q)= D(Q\|P)$, $\hat{D}$ is right-convex.
Hence, it doesn't lose any generality if we only consider the right-convex divergences.

Before we introduce the integral and differential operators, we prove the following important lemma.

\textbf{Lemma 1}.
~\emph{
\label{lem_functional}
 $\int_{0}^t \frac{1}{s} D(P\|R(s)) ds$ and $t\frac{d}{dt}D(P\|R(t))$ don't depend on $t$ explicitly. That is, they only depend on $P$ and $R(t)$.
}

\emph{Proof of Lemma 1.}~
Let $G(P,R(t))\eqdef \int_{0}^1 \frac{1}{s} D(P\|(R(t)-P)s+P) ds=\int_{0}^1 \frac{1}{s} D(P\|R(ts)) ds$.
We prove that $\int_{0}^t \frac{1}{s} D(P\|R(s)) ds=G(P,R(t))$.
When $t=0$, we have
\begin{align}
\int_{0}^t \frac{1}{s} D(P\|R(s)) ds=0=\int_{0}^1 \frac{1}{s} D(P\|P) ds=G(P,R(0)).
\end{align}

When $t>0$, we have 
\begin{align}
G(P,R(t))=\int_{0}^1 \frac{1}{s} D(P\|R(ts)) ds=\int_{0}^{t} \frac{1}{s'} D(P\|R(s')) ds',
\end{align}
where $s'= ts$.
Hence, the result follows.

Next we prove that $t\frac{d}{dt}D(P\|R(t))=H(P,R(t))$, where $H(P,R(t))\eqdef \left.\frac{d}{ds}D(P\|(R(t)-P)s+P)\right|_{s=1}$. 

When $t=0$, we have
\begin{align}
t\frac{d}{dt}D(P\|R(t))=0=\left.\frac{d}{ds}D(P\|P)\right|_{s=1}=H(P,R(0)).
\end{align}

When $t>0$, we have 
\begin{align}
H(P,R(t))=\left.\frac{d}{ds}D(P\|R(ts))\right|_{s=1}=t\left.\frac{d}{d(ts)}D(P\|R(ts))\right|_{s=1}=t\left.\frac{d}{dt'}D(P\|R(t'))\right|_{t'=t},
\end{align}
where $t'= ts$.
Hence, the result follows.

From this lemma, we can introduce integral and differential operators as follows.

\textbf{Definition 1}  (Integral operator).
~\emph{
Let $D$ be a divergence.
When $\frac{1}{t}D(P\|R(t))$ is integrable on $[0,1]$ for the Lebesgue measure, we define an integral operator $\Psi$ as follows.
\begin{align}
\Psi[D](P\|R(t))\eqdef \int_{0}^t \frac{1}{s} D(P\|R(s)) ds. 
\end{align}
}

\textbf{Definition 2} (Differential operator).
~\emph{
Let $D$ be a divergence.
When $D(P\|R(t))$ is differentiable with respect to $t\in[0,1]$, we define a differential operator $\Psi^{-1}$ as follows.
\begin{align}
\Psi^{-1}[D](P\|R(t))\eqdef t\frac{d}{dt}D(P\|R(t)).
\end{align}
}

We can define the integral and differential operators for $D(R(t)\|P)$ in the same way.

From the definitions, we can easily check that $\Psi \circ \Psi^{-1} = \Psi^{-1}\circ \Psi  = 1$.
For $k\geq 0$, we define $\Psi^{k}[D]$ recursively as $\Psi^{k+1}\eqdef \Psi \circ \Psi^{k}$ and $\Psi^{0}$ denotes an identity operator.

\textbf{3. Main Results}

\emph{3.1 Properties of convex divergences}

We show some results for the right-convex divergences.
As we mentioned in the previous section, the same results hold for the left-convex divergences.
In the following, $P,Q\in \mathcal{P}$ denote probability distributions and $R(t)$ denotes mixture distributions $(Q-P)t+P$ for $t\in[0,1]$ as with the previous section.

\textbf{Theorem 1} (Basic theorem). 
~\emph{
Let $t\in[0,1]$.
Let $D$ be a right-convex divergence and let $\frac{1}{t}D(P\|R(t))$ be integrable with respect to $t$.
Then, $\Psi[D](P\|R(t))$ is also a right-convex divergence and 
\begin{align}
\label{ineq_integral_divergence} 
D(P\|R(t)) \geq \Psi[D](P\|R(t))\geq 0.
\end{align}
If $P\neq Q$, the divergence $\Psi[D](P\|R(t))$ is strictly increasing with respect to $t$.
Furthermore, $\frac{1}{t}\Psi[D](P\|R(t))$ is integrable with respect to $t$.
}

\textbf{Theorem 2} (Monotonically decreasing sequence).
~\emph{
Let $t\in[0,1]$.
Let $D$ be a right-convex divergence and let $\frac{1}{t}D(P\|R(t))$ be integrable with respect to $t$.
Then, $\Psi^{k}[D](P\|R(t))$ are right-convex divergences for $k \geq 1$ and $\{\Psi^{k}[D](P\|R(t))\}$ constitute a monotonically decreasing sequence.
\begin{align}
D(P\|R(t)) \geq \Psi[D](P\|R(t)) \geq \Psi^{2}[D](P\|R(t))\geq \cdots \geq  \Psi^{k}[D](P\|R(t)) \cdots \geq 0.
\end{align}
If $P\neq Q$, the divergence $\Psi^{k}[D](P\|R(t))$ are strictly increasing with respect to $t$.}

\textbf{Theorem 3}  (Monotonically decreasing sequence 2).
~\emph{
Let $t\in[0,1]$.
Let $D$ be a differentiable right-convex divergence.
Then, $\Psi^{-1}[D](P\|R(t))$ is a divergence and $\Psi^{k}[D](P\|R(t))$ are right-convex divergence for $k \geq 1$.
They constitute a monotonically decreasing sequence.
\begin{align}
\label{ineq_derivative_divergence}
\Psi^{-1}[D](P\|R(t)) \geq D(P\|R(t)) \geq \Psi[D](P\|R(t)) \geq \Psi^{2}[D](P\|R(t))\geq \cdots \geq  \Psi^{k}[D](P\|R(t)) \cdots \geq 0.
\end{align}
If $P\neq Q$, the divergence $\Psi^{-1}[D](P\|R(t))$ and $\Psi^{k}[D](P\|R(t))$ are strictly increasing with respect to $t$.
}

\textbf{Corollary 1}.
~\emph{
Let $D$ be a differentiable right-convex divergence.
If $P\neq Q$, $D(P\|R(t))$ is strictly increasing with respect to $t$.
}

\emph{Proof.}~
When $P\neq R(t)$, $\Psi^{-1}[D](P\|R(t))=t\frac{d}{dt} D(P\|R(t)) > 0$ holds from Theorem 3.
Since $P\neq R(t)$ holds when $t>0$ and $P\neq Q$, the result follows.

\emph{3.2 Proofs of main results}

We first show the following lemmas.

\textbf{Lemma 2}.
~\emph{
\label{lem_convexity}
Let $x,y\in \mathbb{R}$ and let $g: \mathbb{R}\rightarrow \mathbb{R}$ be a convex function.
Then,
\begin{align}
\label{ineq_convexity}
g(x)\geq g(y)+g'(y)(x-y)
\end{align}
and $g'(x)\leq g'(y)$ for $x<y$.
}

The inequality (\ref{ineq_convexity}) is equivalent to the condition that $g$ is convex (Boyd {\&} Vandenberghe, 2004).
By taking the sum of (\ref{ineq_convexity}) and the inequality with $x$ and $y$ exchanged in (\ref{ineq_convexity}), we have
$(y-x)(g'(y)-g'(x))\geq 0$.

\textbf{Lemma 3}.
~\emph{
\label{lem_convexity2}
If $D$ is a right-convex divergence, $D(P\|R(t))$ is convex in $t$.
}

\emph{Proof.}~
For $\lambda\in[0,1]$, 
\begin{align}
\label{ineq_def_convexity}
\lambda D(P\|R(t_1)) + (1-\lambda)  D(P\|R(t_2)) \geq D(P\|\lambda R(t_1) + (1-\lambda)R(t_2)).
\end{align}
Since $R(t)=(Q-P)t+P$, we have $\lambda R(t_1) + (1-\lambda)R(t_2)=R(\lambda t_1+ (1-\lambda) t_2)$.
By combining this equality and (\ref{ineq_def_convexity}), the result follows.
 \sq

\emph{Proof of Theorem 1.}~
When $P= R(t)$, since $P=R(t)$ holds if and only if $P=Q$ or $t=0$, $\Psi[D](P\|R(t))=0$ holds from the definition of the integral operator $\Psi$.
When $P\neq R(t)$, $\Psi[D](P\|R(t))>0$ holds as well.
Hence, $\Psi[D](P\|R(t))$ is a divergence.

If $P\neq Q$, since $P\neq R(t)$ holds for $t>0$ and $D(P\|R(t))$ is positive, we can easily confirm that $\Psi[D](P\|R(t))$ is strictly increasing with respect to $t$.

We prove the convexity of $\Psi[D](P\|R(t))$.
For $\lambda \in[0,1]$, $R_1(t)=(Q_1-P)t+P$ and $R_2(t)=(Q_2-P)t+P$,
\begin{align}
 \lambda \Psi[D](P\|R_1(t)) + (1-\lambda)\Psi[D](P\|R_2(t))=\int_0^{t} \frac{1}{s} \bigl(\lambda D(P\| R_1(s)) + (1-\lambda)D(P\|R_2(s))\bigr) ds \\ \nonumber
\geq \int_0^{t} \frac{1}{s} D(P\|\lambda R_1(s) + (1-\lambda)R_2(s)) ds =\int_0^{t} \frac{1}{s} D(P\|(\lambda Q_1 + (1-\lambda)Q_2-P)s + P) ds \\ \nonumber
=\Psi[D](P\|(\lambda Q_1 + (1-\lambda)Q_2-P)t + P))=\Psi[D](P\|\lambda R_1(t) + (1-\lambda)R_2(t)),
\end{align}
where we use the convexity of $D$.
Hence, $\Psi[D](P\|R(t))$ is a right-convex divergence.

Next, we prove (\ref{ineq_integral_divergence}).
Let $g(t)\eqdef \Psi[D](P\|R(t))$.

Since $g(t)$ is convex from the convexity of $\Psi[D](P\|R(t))$ and Lemma 3, by applying Lemma 2, we have
\begin{align}
g(0)\geq g(t)+g'(t)(0-t).
\end{align}
Since $g(0)=0$ and $g'(t)=\frac{1}{t}D(P\|R(t))$, we have (\ref{ineq_integral_divergence}). 

Since $\frac{1}{s}D(P\|R(s))$ is integrable, by dividing $D(P\|R(s))\geq \Psi[D](P\|R(s)) \geq 0$ by $s$ and integrating with respect to $s$ from 0 to $t$, we showed that $\frac{1}{s}\Psi[D](P\|R(s))$ is integrable.
\sq

\emph{Proof of Theorem 2.}~
From Theorem 1, the result of applying the integral operator to an integrable right-convex divergence is also an integrable right-convex divergence.
Hence, by applying Theorem 1 repeatedly, the result follows. \sq

\emph{Proof of Theorem 3.}~
We first prove that $\Psi^{-1}[D](P\|R(t)) \geq D(P\|R(t))$.
Let $h(t)\eqdef D(P\|R(t))$.
Since $h(t)$ is convex from Lemma 3, by applying Lemma 2, we have
\begin{align}
h(0)\geq h(t)+h'(t)(0-t).
\end{align}
Since $h(0)=0$ and $h'(t)=\frac{d}{dt}(P\|R(t))$, we have $\Psi^{-1}[D](P\|R(t)) \geq D(P\|R(t))$.

When $P\neq R(t)$, from $\Psi^{-1}[D](P\|R(t)) \geq D(P\|R(t))>0$, we have $\Psi^{-1}[D](P\|R(t))>0$.
Since $P=R(t)$ holds if and only if $t=0$ or $P=Q$, $t\frac{d}{dt} D (P\|R(t)) =0$ for $t=0$ and $t\frac{d}{dt} D (P\|P) =0$ for $P=Q$.
Hence, $\Psi^{-1}[D](P\|R(t))$ is a divergence.

For $t_2>t_1\geq 0$, 
\begin{align}
\Psi^{-1}[D](P\|R(t_2))-\Psi^{-1}[D](P\|R(t_1))=(t_2-t_1)h'(t_2) + t_1(h'(t_2)-h'(t_1)).
\end{align}
When $P\neq R(t)$, the first term in RHS is positive from $h'(t)=\frac{1}{t}\Psi^{-1}[D](P\|R(t))$.
From Lemma 2, because the second term in RHS is non-negative, $\Psi^{-1}[D](P\|R(t))$ is strictly increasing.

Since $D(P\|R(t))$ is continuous from the differentiable assumption and $\lim_{\epsilon\rightarrow +0} \frac{D(P\|R(\epsilon))}{\epsilon}=\frac{d}{dt}D(P\|R(t))|_{t=0} < \infty$, $\frac{1}{t}D(P\|R(t))$ is integrable.

By applying Theorem 2, we can prove the rest part of the theorem. \sq

\textbf{4. Examples of divergence sequences}

In this section, we show a specified example of divergence sequences by applying theorems in the previous section.

\emph{4.1 Convex divergence sequences that include the KL-divergence} 

We introduce new divergences by using the polylogarithm and we show that they constitute monotonically decreasing sequence that include the KL-divergence or the reverse KL-divergence.

\textbf{Definition 3}.~
\emph{
We define $\mathrm{PL}_k: \mathcal{P}\times \mathcal{P}\rightarrow\mathbb{R}$ for $k\geq 0$ as follows.
\begin{align}
\mathrm{PL}_k(P\|Q)\eqdef \int p\mathrm{Li}_k\biggl(1-\frac{q}{p}\biggr) d\mu,
\end{align}
where $\mathrm{Li}_k(z)$ is the polylogarithm.
The polylogarithm is defined as $\mathrm{Li}_{0}(z)\eqdef \frac{z}{1-z}$ and 
\begin{align}
\label{eq_def_polylog_integral}
\mathrm{Li}_k(z)\eqdef \frac{z}{\Gamma(k)}\int_0^\infty \frac{x^{k-1}}{e^x-z}  dx,
\end{align}
for $k>0$, where $z$ is the complex argument.
}

The polylogarithm satisfies 
\begin{align}
\label{eq_def_polylog}
\mathrm{Li}_{k+1}(z)=\int_0^z \frac{\mathrm{Li}_k(x)}{x}dx, 
\end{align}
\begin{align}
\label{eq_poly_0}
\mathrm{Li}_k(0)=0
\end{align}
and 
\begin{align}
\mathrm{Li}_1(z)=-\log(1-z).
\end{align}

We name $\mathrm{PL}_k$ after ``PolyLogarithm''.

\textbf{Proposition 1}.~
\emph{
For $k\geq 0$, $\mathrm{PL}_k(P\|R(t))=\Psi^k[\mathrm{PL}_0](P\|R(t))$.
Furthermore, $\mathrm{PL}_k(P\|R(t))$ are right-convex divergences and $\{\mathrm{PL}_k(P\|R(t))\}$ constitute a monotonically decreasing sequence.
}

From the definition of the polylogarithm and $\mathrm{PL}_k$, we can easily confirm that 
\begin{align}
\label{eq_chi_square}
\mathrm{PL}_0(P\|R(t))=\int \frac{p(p-r(t))}{r(t)} d\mu=\int \frac{(r(t)-p)^2}{r(t)} d\mu=\chi^2(P\|R(t))
\end{align} and 
\begin{align}
\mathrm{PL}_1(P\|R(t))=\mathrm{KL}(P\|R(t)),
\end{align}
where $\chi^2(P\|Q)\eqdef \int \frac{(q-p)^2}{q} d\mu$ is the Neyman $\chi^2$-divergence.
We have derived a lower bound for the KL-divergence by using the result of this theorem $\mathrm{KL}(P\|R(t))=\Psi[\chi^2](P\|R(t))$ in (Nishiyama, 2019).

\textbf{Definition 4}.~
\emph{
We define $\mathrm{SL}_k: \mathcal{P} \times \mathcal{P}\rightarrow\mathbb{R}$ for $k\geq 0$ as follows.
For $k=0$, 
\begin{align}
\mathrm{SL}_0(P\|Q)\eqdef \mathrm{J}(P,Q).
\end{align}
For $k\geq 1$,
\begin{align}
\label{def_eq_sl}
\mathrm{SL}_k(P\|Q)\eqdef \mathrm{J}(P,Q) - \sum_{j=1}^k \mathrm{PL}_j(P\|Q) ,
\end{align}
where $\mathrm{J}(P,Q)$ is the Jeffreys divergence defined as $\mathrm{J}(P,Q)\eqdef \mathrm{KL}(P\|Q) + \mathrm{KL}(Q\|P)$ (Jeffreys, 1946).
}

We name $\mathrm{SL}_k$ after ``Sum of polyLogarithm''.

Since $\mathrm{PL}_1(P\|R(t))=\mathrm{KL}(P\|R(t))$, we have $\mathrm{SL}_1(P\|R(t))=\mathrm{KL}(R(t)\|P)$.

\textbf{Proposition 2}.~
\emph{
For $k\geq 0$, $\mathrm{SL}_k(P\|R(t))=\Psi^k[\mathrm{SL}_0](P\|R(t))$.
For $k\geq 1$, $\mathrm{SL}_k(P\|R(t))$ are right-convex divergences and $\{\mathrm{SL}_k(P\|R(t))\}$ constitute a monotonically decreasing sequence for $k\geq 0$.
}

Hence, the right-convex divergence sequences $\{\mathrm{PL}_k\}$ and $\{\mathrm{SL}_k\}$ include the KL-divergence and the reverse KL-divergence, respectively.

\emph{4.2 Squared Hellinger distance}

The squared Hellinger distance is defined as $\mathrm{Hel}^2(P,Q)\eqdef \frac{1}{2} \int (\sqrt{q}-\sqrt{p})^2 d\mu$.
Since the squared Hellinger distance belongs to the $f$-divergence and differentiable, it  satisfies assumptions in Theorem 3.
We calculate $\Psi^{-1}[\mathrm{Hel}^2]$ and $\Psi[\mathrm{Hel}^2]$.

 $\Psi^{-1}[\mathrm{Hel}^2]$ is 
\begin{align}
\Psi^{-1}[\mathrm{Hel}^2](P\|R(t))&=\frac{1}{2}\int (r(t)-p)\frac{\sqrt{r(t)}-\sqrt{p}}{\sqrt{r(t)}} d\mu=\frac{1}{2}\int (\sqrt{r(t)}+\sqrt{p})\frac{(\sqrt{r(t)}-\sqrt{p})^2}{\sqrt{r(t)}} d\mu \\ \nonumber
&= \mathrm{Hel}^2(P,R(t))+\frac{1}{2}\int (\sqrt{r(t)}-\sqrt{p})^2\sqrt{\frac{p}{r(t)}} d\mu.
\end{align}

$\Psi[\mathrm{Hel}^2]$ is 
\begin{align}
\Psi[\mathrm{Hel}^2](P\|R(t))&=\frac{1}{2}\int\int_0^t \frac{1}{s}(\sqrt{r(s)}-\sqrt{p})^2 dsd\mu =\frac{1}{2}\int\int_p^{r(t)} \frac{1}{s'-p}(\sqrt{s'}-\sqrt{p})^2 ds'd\mu \\ \nonumber
&=\frac{1}{2}\int\int_p^{r(t)} \frac{\sqrt{s'}-\sqrt{p}}{\sqrt{s'}+\sqrt{p}} ds'd\mu,
\end{align}
where $s'=r(s)=(q-p)s+p$.

Since $\int \frac{\sqrt{x}-\sqrt{p}}{\sqrt{x}+\sqrt{p}} dx= -4\sqrt{p}\sqrt{x} +4p\log(\sqrt{x}+\sqrt{p})+x+\mathrm{const.}$, we have
\begin{align}
\Psi[\mathrm{Hel}^2](P\|R(t))&=2-2\int\sqrt{p}\sqrt{r(t)} d\mu+2\int p\log\frac{(\sqrt{p}+\sqrt{r(t)})}{2\sqrt{p}}d\mu \\ \nonumber
&=2\mathrm{Hel}^2(P\|R(t))+2\int p\log\frac{(\sqrt{p}+\sqrt{r(t)})}{2\sqrt{p}}d\mu,
\end{align}  
where we use $\int pd\mu=\int r(t)d\mu=1$.

\emph{4.3 Proofs of propositions}

\emph{Proof of Proposition 1.}~ 
First, we show that $\frac{1}{t}\mathrm{PL}_0(P\|R(t))$ is integrable.
From (\ref{eq_chi_square}), since $\mathrm{PL}_0(P\|R(t))=\int \frac{(r(t)-p)^2}{r(t)} d\mu = \chi^2(P\|R(t))$ holds,
\begin{align}
\int_0^t \frac{1}{s}\mathrm{PL}_0(P\|R(s))ds&=\int_0^t \frac{1}{s} \int \frac{(r(s)-p)^2}{r(s)} d\mu ds=\int (q-p)\int_0^t \frac{r(s)-p}{r(s)}ds d\mu \\  \nonumber
 &=\int \int_p^{r(t)} \frac{s'-p}{s'}ds' d\mu = \int p\log\frac{p}{r(t)}d\mu = \mathrm{KL}(P\|R(t)), 
\end{align}
where $s'=r(s)=(q-p)s+p$ and we use $\int p d\mu=\int r(t) d\mu=1$.

In addition, since the $\chi^2$-divergence belongs to the $f$-divergence,  $\mathrm{PL}_0(P\|R(t))$ satisfies assumptions in Theorem 2.

Hence, if we prove that $\mathrm{PL}_k(P\|R(t))=\Psi^k[\mathrm{PL}_0](P\|R(t))$, we can also show the rest part of the proposition from Theorem 2.

We prove this equality by the induction.

The case $k=0$ is trivial.

Suppose that $\mathrm{PL}_l(P\|R(t))=\Psi^l[\mathrm{PL}_0](P\|R(t))$ for $k=l$, we prove the same equality for $k=l+1$.
From the assumption of the induction, 
\begin{align}
\label{eq_induction1}
\Psi^{l+1}[\mathrm{PL}_0](P\|R(t))&=\Psi[\mathrm{PL}_l](P\|R(t))=\int_0^t \frac{1}{s}\int p\mathrm{Li}_l\biggl(1-\frac{r(s)}{p}\biggr) d\mu ds =\int p\int_0^t \frac{1}{s}\mathrm{Li}_l\biggl(1-\frac{r(s)}{p}\biggr) ds d\mu.
\end{align}
When $p\neq r(s)$, by putting $s'=1-\frac{r(s)}{p}=-\frac{(q-p)s}{p}$ and using (\ref{eq_def_polylog}), we have

\begin{align}
\int_0^t \frac{1}{s}\mathrm{Li}_l\biggl(1-\frac{r(s)}{p}\biggr) ds=\int_0^{1-\frac{r(t)}{p}} \frac{1}{s'} \mathrm{Li}_l(s') ds' =\mathrm{Li}_{l+1} \biggl(1-\frac{r(t)}{p}\biggr).
\end{align}
When $p=r(s)$, from (\ref{eq_poly_0}), the same equality holds. 
By substituting this equality into (\ref{eq_induction1}), we have
\begin{align}
\Psi^{l+1}[\mathrm{PL}_0](P\|R(t))=\int p\mathrm{Li}_{l+1} \biggl(1-\frac{r(t)}{p}\biggr)d\mu=\mathrm{PL}_{l+1}(P\|R(t)).
\end{align}

Then, the result follows.   \sq

\emph{Proof of Proposition 2.}~ 
As mentioned in the previous subsection, $\mathrm{SL}_1(P\|R(t))=\mathrm{KL}(R(t)\|P)$ holds.
Since the reverse KL-divergence is differentiable and belongs to the $f$-divergence, $\mathrm{SL}_1(P\|R(t))$ satisfies assumptions in Theorem 3.   

If we prove that $\mathrm{SL}_1(P\|R(t))=\Psi[\mathrm{SL}_0](P\|R(t))$ and $\mathrm{SL}_k(P\|R(t))=\Psi^{k-1}[\mathrm{SL}_1](P\|R(t))$ for $k\geq 1$, we can prove the whole proposition from Theorem 3.

First, we show that $\mathrm{SL}_1(P\|R(t))=\Psi[\mathrm{SL}_0](P\|R(t))$.
Since the Jeffreys divergence is written as $J(P,Q)=\mathrm{KL}(P\|Q)+\mathrm{KL}(Q\|P)=\int (p-q)\log \frac{p}{q} d\mu$,
\begin{align}
\Psi[\mathrm{SL}_0](P\|R(t))&=\int_0^t \int \frac{1}{s} (r(s)-p)\log\frac{r(s)}{p} d\mu ds=\int(q-p)\int_0^t \log\frac{r(s)}{p} dsd\mu \\ \nonumber
&=\int p\int_1^{\frac{r(t)}{p}} \log s' ds'd\mu=\int r(t)\log\frac{r(t)}{p} d\mu=\mathrm{KL}(R(t)\|P)=\mathrm{SL}_1(P\|R(t)),
\end{align}
where $s'=\frac{r(s)}{p}=\frac{(q-p)s}{p}+1$ and we use $\int pd\mu=\int r(t) d\mu=1$.
Hence, we have the result.
From the definition of $\mathrm{SL}_1$, we have 
\begin{align}
\label{eq_psi_J}
\Psi[\mathrm{J}](P\|R(t))=\Psi[\mathrm{SL}_0](P\|R(t))=\mathrm{SL}_1(P\|R(t))=\mathrm{J}(P,R(t))-\mathrm{PL}_1(P\|R(t)).
\end{align}

Next, we prove $\mathrm{SL}_k(P\|R(t))=\Psi^{k-1}[\mathrm{SL}_1](P\|R(t))$ for $k\geq 1$ by the induction.

The case $k=1$ is trivial.
Suppose $\mathrm{SL}_l(P\|R(t))=\Psi^{l-1}[\mathrm{SL}_1](P\|R(t))$ for $k=l\geq 1$, we prove the same equality for $k=l+1$.
From Proposition 1, we have $\Psi[\mathrm{PL}_k](P\|R(t))=\mathrm{PL}_{k+1}(P\|R(t))$.
By combining this equality, (\ref{def_eq_sl}) and (\ref{eq_psi_J}), we have
\begin{align}
\Psi^{l}[\mathrm{SL}_1](P\|R(t))&=\Psi[\mathrm{SL}_l](P\|R(t))=\mathrm{J}(R(t), P)-\mathrm{PL}_1(P\|R(t)) - \sum_{j=1}^{l} \mathrm{PL}_{j+1}(P\|R(t))=\mathrm{SL}_{l+1}(P\|R(t)).
\end{align}
Then, the result follows.    \sq

\textbf{5. Conclusion}

We focused on the convex divergences and showed some properties through the integral and differential operators.
The integral operator preserves the convexity and the properties of divergence, and the differential operator preserves the properties of divergence.
We showed that the convex divergence applied the integral operator multiple times constitute a monotonically decreasing sequence of convex divergences.
In addition, we defined new decreasing sequences of divergences that include the KL-divergence, the reverse KL-divergence, the $\chi^2$-divergence and the Jeffreys divergence by using the polylogarithm.

It is our future work to study properties of the convex divergences in more detail.


\textbf{References}

{\leftskip=2em\parindent=-2em

Ajjanagadde, G., Makur, A., Klusowski, J., {\&} Xu, S. (2017). Lecture notes on information theory.

Amari, S. I., {\&}  Cichocki, A. (2010). Information geometry of divergence functions. \emph{Bulletin of the Polish Academy of Sciences: Technical Sciences, 58}(1), 183-195.
https://doi.org/10.2478/v10175-010-0019-1

Billingsley, P. (2008). \emph{Probability and measure}. John Wiley {\&} Sons.

Boyd, S., {\&} Vandenberghe, L. (2004). \emph{Convex optimization}. Cambridge university press.

Bregman, L. M. (1967). The relaxation method of finding the common point of convex sets and its application to the solution of problems in convex programming. \emph{USSR computational mathematics and mathematical physics, 7}(3), 200-217.
https://doi.org/10.1016/0041-5553(67)90040-7

Cover, T. M., {\&} Thomas, J. A. (2012). \emph{Elements of information theory}. John Wiley {\&} Sons.

Csisz{\'a}r, I., {\&} Shields, P. C. (2004). Information theory and statistics: A tutorial. \emph{Foundations and $Trends^{TM}$ in Communications and Information Theory, 1}(4), 417-528.
https://doi.org/10.1561/0100000004 

Jeffreys, H. (1946). An invariant form for the prior probability in estimation problems. \emph{Proceedings of the Royal Society of London. Series A. Mathematical and Physical Sciences, 186}(1007), 453-461.
https://doi.org/10.1098/rspa.1946.0056

Kullback, S. (1997). \emph{Information theory and statistics}. Courier Corporation.

Nishiyama, T. (2018). Divergence Network: Graphical calculation method of divergence functions. \emph{arXiv preprint arXiv:1810.12794}.

Nishiyama, T. (2019). A New Lower Bound for Kullback-Leibler Divergence Based on Hammersley-Chapman-Robbins Bound. \emph{arXiv preprint arXiv:1907.00288}.

Sason, I., {\&} Verdu, S. (2016). $ f $-divergence Inequalities. \emph{IEEE Transactions on Information Theory, 62}(11), 5973-6006.
https://doi.org/10.1109/tit.2016.2603151

Van Erven, T., {\&} Harremos, P. (2014). R{\'e}nyi divergence and Kullback-Leibler divergence. \emph{IEEE Transactions on Information Theory, 60}(7), 3797-3820.
https://doi.org/10.1109/tit.2014.2320500

Wood, D. C. (1992). The computation of polylogarithms.

\end{document}